\documentclass{amsart}
\usepackage{amssymb,euscript,amsmath, mathrsfs}
\usepackage[dvips]{graphicx}
\usepackage[dvips]{color}

\newcounter{ENUM}

\input xy
\xyoption{all}
\CompileMatrices

\def\<{\langle}
\def\>{\rangle}
\def\0{{{\bf 0}}}

\def\OO{{\mathcal O}}

\def\CK{{\mathcal K}}

\def\CP{{\mathcal P}}

\def\tkappa{{\tilde \kappa}}

\def\CC{{\mathbb C}}

\def\FF{{\mathbb F}}

\def\QQ{{\mathbb Q}}

\def\tx{{\tilde x}}

\def\tpi{{\tilde \pi}}

\def\trho{{\tilde \rho}}

\def\rhobar{{\overline{\rho}}}
\def\rhobar{{\overline{\rho}}}
\def\pibar{{\overline{\pi}}}

\def\unif{{\varpi}}

\newcommand{\cInd}{\operatorname{c-Ind}}

\newcommand{\St}{\operatorname{St}}

\newcommand{\LL}{\operatorname{LL}}
\newcommand{\Rep}{\operatorname{Rep}}

\newcommand{\soc}{\operatorname{soc}}

\def\Hom{\operatorname{Hom}}

\def\End{\operatorname{End}}

\def\GL{\operatorname{GL}}

\def\Spec{\operatorname{Spec}}

\newcommand{\margh}[1]{}

\newtheorem{thm}{Theorem}[section]
\newtheorem{theorem}[thm]{Theorem}

\newtheorem{proposition}[thm]{Proposition}
\newtheorem{lemma}[thm]{Lemma}

\newtheorem{corollary}[thm]{Corollary}
\newtheorem{conj}[thm]{Conjecture}

\theoremstyle{definition}

\newtheorem{definition}[thm]{Definition}

\newtheorem{example}[thm]{Example}

\newtheorem{remark}[thm]{Remark}

\numberwithin{equation}{section}

\begin{document}
\title{Whittaker models and the integral Bernstein center for $\GL_n$}
\author{David Helm}

\begin{abstract}
We establish integral analogues of results of Bushnell and Henniart~\cite{BH-whittaker}
for spaces of Whittaker functions arising from the groups $\GL_n(F)$ for $F$
a $p$-adic field.  We apply the resulting theory to the existence of representations
arising from the conjectural ``local Langlands correspondence in families'' of~\cite{emerton-helm},
and reduce the question of the existence of such representations to a natural
conjecture relating the integral Bernstein center of $\GL_n(F)$ to the
deformation theory of Galois representations.
\end{abstract}

\maketitle

\section{Introduction}

The {\em Bernstein center} is a central tool in the study of the smooth
complex representations of a $p$-adic group $G$.  Introduced by Bernstein
and Deligne in~\cite{BD}, the Bernstein center is a commutative ring that
acts naturally on every smooth complex representation of $G$.  Its primitive
idempotents thus decompose the category $\Rep_{\CC}(G)$ of smooth complex
representations of $G$ into full subcategories
known as ``blocks''; any object of $\Rep_{\CC}(G)$ then has a canonical
decomposition as a product of factors, one from each block.  Moreover, Bernstein
and Deligne give a description of each block, and show that the action of
the Bernstein center on each block factors through a finite type $\CC$-algebra
(the center of the block), that can be described in a completely explicit fashion.

One advantage of this approach is that it allows one to give purely algebraic
proofs of results that were classically proven using difficult technqiues from
Harmonic analysis.  For example, one can regard a (not necessarily admissible) smooth
representation of $G$ as a sheaf on the spectrum of the Bernstein center; doing
so provides a purely algebraic analogue of the Fourier decomposition of this representation
as a ``direct integral'' of irreducible representations over a suitable measure space.
A clear benefit of this algebraic approach is that it applies even when one
considers representations over fields (or even rings) other than the complex numbers.

In~\cite{BH-whittaker}, Bushnell and Henniart applied the above ideas to the
study of Whittaker models in the complex representation theory of $p$-adic groups.
In particular, they study the space $\cInd_U^G \Psi$, where $U$ is the unipotent
radical of a Borel subgroup $B$ of $G$, and $\Psi$ is a ``generic character'' of $U$.
This space is dual to the space of functions in which Whittaker models live.  They establish
two key technical results about this space: first, if $e$ is an idempotent of the
Bernstein center corresponding to a particular block of $\Rep_{\CC}(G)$,
they show that $e \cInd_U^G \Psi$ is finitely generated as a $\CC[G]$-module.
Second, they show that the center of the block corresponding to $e$ acts faithfully
on $e \cInd_U^G \Psi$, and that in many clases (in particular $G = \GL_n$), this
center is the full endomorphism ring of $e \cInd_U^G \Psi$.  As an application,
they give a purely algebraic proof of a vanishing theorem, originally due to
Jacquet, Piatetski-Shapiro, and Shalika (\cite{JPS}, Lemma 3.5), for 
functions in $\cInd_U^G \Psi$.

The theory of the Bernstein center for categories of representations
over coefficient rings other than $\CC$ presents significant technical difficulties
beyond the complex case; in particular the approach of Bernstein and Deligne
makes heavy use of certain properties of cuspidal representations that only hold
over fields of characteristic zero.  Dat~\cite{dat-integral} has established some basic
structural results for general $p$-adic groups in this context, but for general
groups not much more is known.

For $\GL_n(F)$, however, previous work of the author~\cite{bernstein1} establishes
a much more detailed structure theory for the center of the category $\Rep_{W(k)}(\GL_n(F))$
of smooth $W(k)[\GL_n(F)]$-modules, where $k$ is an algebraically closed field
of characteristic $\ell$.  Although it does not seem possible to give as complete and
thorough description of the center in this context, the results of~\cite{bernstein1}
are detailed and precise enough that one might hope for applications in both
representation theory and arithmetic.  In particular it is natural to ask whether
one can apply the theory of~\cite{bernstein1} to analogues of the questions
studied in~\cite{BH-whittaker}.

The purpose of this paper is twofold: first, to establish integral versions
of the results of~\cite{BH-whittaker} discussed above in the case where $G = \GL_n(F)$,
using the results and techniques of~\cite{bernstein1}.  The second is to apply
these techniques to questions concerning the local Langlands correspondence of~\cite{emerton-helm}.

Rather than attempting to mimic the arguments of~\cite{BH-whittaker}, our
approach to the first goal is by a quite different argument.  We make
heavy use of the fact that the results of~\cite{BH-whittaker} hold over $\CC$,
and therefore also over the field $\overline{\CK}$, where $\CK$ is the field
of fractions of $W(k)$.  Our approach here relies on the Bernstein-Zelevinski
theory of the derivative, developed over $W(k)$ in~\cite{emerton-helm}, and in
particular on computing the derivatives of certain projective objects of
$\Rep_{W(k)}(G)$ first considered in~\cite{bernstein1}.  This computation makes
heavy use of the notion of an essentially AIG representation, introduced
in~\cite{emerton-helm}.

Our main results about $\cInd_U^G \Psi$ are established in section~\ref{sec:admissible}.
As an application we prove a basic fact about essentially AIG representations
that was conjectured in~\cite{emerton-helm}.

We now discuss our second goal in more detail.  In~\cite{emerton-helm}, Emerton and the
author introduce a conjectural ``local Langlands correspondence in families''.  The main
result states roughly that given a Galois representation $\rho: G_F \rightarrow \GL_n(A)$,
for $A$ a suitable complete local $W(k)$-algebra, there is at most one admissible
$A[\GL_n(F)]$-module $\pi(\rho)$, satisfying a short list of technical conditions, such
that at characteristic zero points ${\mathfrak p}$ of $\Spec A$, the representations
$\pi(\rho)_x$ and $\rho_x$ are related by a variant of the local Langlands correspondence.
(We refer the reader to Theorem~\ref{thm:E-H} for a precise statement).

Emerton has shown~\cite{emerton-lg} that certain spaces that arise by
considering the completed cohomology of the tower for modular curves have a natural
tensor factorization, and that the tensor factors that arise in this way are isomorphic
to $\pi(\rho)$ for certain representations $\rho$ of $G_{\QQ_p}$ over Hecke algebras.
This result is crucial to his approach to the Fontaine-Mazur conjecture.

However, the results of~\cite{emerton-helm} do not address the question of whether
the representations $\pi(\rho)$ exist in general.  They also leave several fundamental
questions about the structure of $\pi(\rho)$ unanswered.

Our approach to these questions revolves around the concept of a
``co-Whittaker'' module, introduced in section~\ref{sec:whittaker}.
The key point is that the families $\pi(\rho)$, when they exist, are
co-Whittaker modules.  Moreover, there is a natural connection
between co-Whittaker modules and the module $\cInd_U^G \Psi$, which can
be regarded as a ``universal'' co-Whittaker module, in a sense we make
precise in Theorem~\ref{thm:universal}.

In the final section we apply our structure theory to the local Langlands
correspondence in families.  In particular, we reduce the construction
of the families $\pi(\rho)$, where $\rho$ is a deformation of a given
Galois representation $\rhobar$ over $k$, to the question of the
existence of a certain map from the integral Bernstein center to
the universal framed deformation ring of $\rhobar$ (Theorem~\ref{thm:main}).
Moreover, we show that the converse also holds (conditionally on a result
from the forthcoming work~\cite{bernstein2}.)

There are several advantages to reformulating the question of the existence
of $\pi(\rho)$ in this way.  The first is that the new conjectures that
arise in this way (Conjecture~\ref{conj:weak} and the stronger variant
Conjecture~\ref{conj:strong}) are interesting in their own right.  
Indeed, one can regard Conjecture~\ref{conj:strong} as a geometric
reformulation of the correspondence between supercuspidal support
(for admissible representations of $G$) and semisimplification
of Galois representations (see remark~\ref{rem:interpolation}.)

Moreover, it seems to be substantially easier in practice
to construct maps from the Bernstein center to universal framed
deformation rings than it is to directly construct representations
$\pi(\rho)$.  For example, it is not difficult to show
that both conjectures hold after inverting $\ell$, as well as in
the case where $\ell$ is a {\em banal} prime; that is, when
the order $q$ of the residue field of $F$ has the property that
$1, q, \dots, q^n$ are distinct mod $\ell$.  (We sketch a proof of
this in Example~\ref{ex:banal}; the details, as well
as deeper results about the conjectures in question, will appear
in the forthcoming work~\cite{bernstein2}.)  In particular, this
approach provides the only currently known construction of $\pi(\rho)$
in the banal setting for $n > 2$.

For small $n$, it is in principle possible to construct suitable
families $\pi(\rho)$ by ad-hoc methods.  It seems unlikely, however,
that these methods can be pushed much beyond the case $n=2$.
For instance, when $n = 2$, $\ell$ is odd, and $q$ is congruent
to $-1$ modulo $\ell$, it is possible to explicitly construct
the representation $\pi(\rho)$, where $\rho$ is the universal
framed deformation ring of $1 \oplus \omega$ (here $\omega$
is the mod $\ell$ cyclotomic character).  This construction
requires a delicate analysis of congruences between lattices
in cuspidal and Steinberg representations of $\GL_2(F)$.
By contrast, if one knows Conjecture~\ref{conj:weak} in this setting,
the construction is clean and straightforward.  In a sense,
the map in Conjecture~\ref{conj:weak} encodes all of these congruences
in a sufficiently systematic way that one does not need to work
with them directly.

The author is grateful to Matthew Emerton, Richard Taylor, and Sug-Woo
Shin for their continued interest and helpful conversation on the ideas
presented here.  This research was partially supported by NSF grant
DMS-1161582.

\section{The integral Bernstein center}

We now fix notation and summarize some of the basic properties of the
integral Bernstein center from~\cite{bernstein1} that will be used throughout
the paper.  Let $p$ and $\ell$ be distinct primes, and let $F$ be a finite
extension of $\QQ_p$.  We let $G$ denote the group $\GL_n(\QQ_p)$.

Let $k$ be an algebraically closed field of characteristic $\ell$,
and let $\CK$ be the field of fractions of $W(k)$.  We will be concerned
with the categories $\Rep_k(G)$, $\Rep_{\CK}(G)$, and
$\Rep_{W(k)}(G)$ of smooth $k[G]$-modules, smooth $\CK[G]$-modules,
and smooth $W(k)[G]$-modules, respectively.

The blocks of $\Rep_{W(k)}[G]$ are parameterized by equivalence classes
of pairs $(L,\pi)$, where $L$ is a Levi subgroup of $G$ and $\pi$ is an
irreducible supercuspidal representation of $L$ over $k$.  Two pairs
$(L,\pi)$ and $(L',\pi')$ are said to be {\em inertially equivalent}
if there is an element $g$ of $G$ such that $L' = g L g^{-1}$ and $\pi'$
is a twist of $\pi^g$ by an unramified character of $L'$.

Given a pair $[L,\pi]$, up to inertial equivalence, we can consider the
full subcategory $\Rep_{W(k)}(G)_{[L,\pi]}$ of $\Rep_{W(k)}(G)$ consisting of
smooth $W(k)[G]$-modules $\Pi$ such that every simple subquotient of $\Pi$
has {\em mod $\ell$ inertial supercuspidal support} (in the sense of~\cite{bernstein1},
Definition 4.10) given by the pair $(L,\pi)$.  By~\cite{bernstein1}, Theorem 10.8,
$\Rep_{W(k)}(G)_{[L,\pi]}$ is a block of $\Rep_{W(k)}(G)$.  For any
smooth $W(k)[G]$-module $\Pi$, we may thus speak of the factor $\Pi_{[L,\pi]}$
of $\Pi$ that lies in $\Rep_{W(k)}(G)_{[L,\pi]}$.

We may also consider the center $A_{[L,\pi]}$ of $\Rep_{W(k)}(G)_{[L,\pi]}$.
We have the following basic structure theory of $A_{[L,\pi]}$:

\begin{theorem}[\cite{bernstein1}, Theorem 12.1]
The ring $A_{[L,\pi]}$ is a finitely generated, reduced, $\ell$-torsion free
$W(k)$-algebra.
\end{theorem}
 
Let $\kappa$ be a $W(k)$-algebra that is a field; then any smooth representation
$\Pi$ of $G$ over $\kappa$ lies in $\Rep_{W(k)}(G)$.  If $\Pi$ is absolutely irreducible
then it lies in $A_{[L,\pi]}$ for some $(L,\pi)$; by Schur's Lemma the action of
$A_{[L,\pi]}$ on $\Pi$ is via a homomorphism $f_{\Pi}: A_{[L,\pi]} \rightarrow \kappa$.

\begin{theorem} \label{thm:bernstein action}
Let $\Pi_1$ and $\Pi_2$ be two absolutely irreducible representations of $G$ over $\kappa$
that lie in $\Rep_{W(k)}(G)_{[L,\pi]}$,
and let $f_1$ and $f_2$ be the maps: $A_{[L,\pi]} \rightarrow \kappa$ giving the action
of $A_{[L,\pi]}$ on $\Pi_1$ and $\Pi_2$ respectively.  Then $f_1 = f_2$ if, and only if,
$\Pi_1$ and $\Pi_2$ have the same supercuspidal support.
\end{theorem}
\begin{proof}
When $\kappa$ has characteristic zero this follows from the classical theory of Bernstein and
Deligne.  When $\kappa$ has characteristic $\ell$ this is (a slight generalization of)
Theorem 12.2 of~\cite{bernstein1}; the proof of that theorem works in
the generality claimed here.
\end{proof}

\section{Essentially AIG representations} \label{sec:derivative}

Our attempts to generalize results from~\cite{BH-whittaker} will rely
on a version of the Bernstein-Zelevinski
theory of the derivative and its related functors that makes sense
for integral representations.
We refer the reader to~\cite{emerton-helm}, section 3.1, for the details of 
this theory over $W(k)$.  Here we will content ourselves with the theory
of the ``top derivative''.

Let $U$ be the unipotent radical of a Borel subgroup of $G$, and let
$\Psi: U \rightarrow W(k)^{\times}$ be a generic character.  Let $W$ be
the module $\cInd_U^G \Psi$; $W$ is independent, up to isomorphism, of the generic character
$\Psi$.  As $U$ has order prime to $\ell$ and $\cInd$ takes projectives to projectives,
$W$ is a projective $W(k)[G]$-module.  The Bernstein factor
$W_{[L,\pi]}$ is then a projective object of $\Rep_{W(k)}(G)_{[L,\pi]}$,
which is closely related to the theory of Whittaker models.

The module $W_{[L,\pi]}$ will be crucial for our
approach to the conjectures of~\cite{emerton-helm} on the local Langlands
correspondence for families of Galois representations.   First, however, we 
need some basic results on the structure of $W_{[L,\pi]}$.
If we invert $\ell$, this was studied systematically by Bushnell-Henniart.  From
our perspective a key result of theirs is:

\begin{theorem}
The natural map 
$$A_{[L,\pi]} \otimes_{W(k)} \CK \rightarrow \End_{\CK[G]}(W_{[L,\pi]} \otimes_{W(k)} \CK)$$
is an isomorphism.  In particular, $\End_{W(k)[G]}(W_{[L,\pi]})$ is commutative.
\end{theorem}
\begin{proof}
We can check the first statement after base change from $\CK$ to $\overline{\CK}$; it then
follows from~\cite{BH-whittaker}, Theorem 4.3.  The second statement is immediate, as
$W_{[L,\pi]}$ is $\ell$-torsion free and thus its endomorphism ring embeds in
$\End_{\CK[G]}(W_{[L,\pi]} \otimes_{W(k)} \CK)$.
\end{proof}

To go further,
we recall that we have an ``$n$th derivative functor'' $V \mapsto V^{(n)}$ from $W(k)[G]$-modules
to $W(k)$ modules with the following properties:

\begin{enumerate}
\item If $V$ is an irreducible $k[G]$-module, then $V^{(n)}$ is a one-dimensional $k$-vector
space if $V$ is generic, and zero otherwise.
\item There is a natural isomorphism of functors $\Hom_{W(k)[G]}(W,-) \rightarrow (-)^{(n)}$.
\item The functor $V \mapsto V^{(n)}$ is exact.
\item If $V$ is an $A[G]$-module for some $W(k)$-algebra $A$, and $B$ is an $A$-algebra,
$(V \otimes_A B)^{(n)}$ is naturally isomorphic to $V^{(n)} \otimes_A B$.
\item If $V$ and $W$ are $A[\GL_n(F)]$ and $A[\GL_m(F)]$-modules, respectively, then there is a natural
isomorphism:
$$[i_P^{\GL_{n+m}(F)} V \otimes W]^{(n+m)} \cong V^{(n)} \otimes W^{(m)}.$$
\item If $V$ is an $A[\GL_n(F)]$-module, there is a natural $A$-linear map $V \rightarrow V^{(n)}$.
\end{enumerate}

\begin{lemma} \label{lemma:freeness}
The top derivative $W^{(n)}_{[L,\pi]}$ is free of rank one over $\End_{W(k)[G]}(W_{[L,\pi]})$.
\end{lemma}
\begin{proof}
We have a natural isomorphism:
$$W^{(n)}_{[L,\pi]} \rightarrow \Hom_{W(k)[G]}(W,W_{[L,\pi]}),$$
and the latter is clearly isomorphic to $\End_{W(k)[G]}(W_{[L,\pi]})$.
\end{proof}

Our goal is to apply this structure theory to ideas from~\cite{emerton-helm}.
We first recall the definition of an essentially AIG representation~\cite{emerton-helm}, 3.2.1.
\begin{definition}
Let $\kappa$ be a $W(k)$-algebra that is a field.  A $\kappa[G]$-module $V$ is {\em essentially AIG} if:
\begin{enumerate}
\item the socle of $V$ is absolutely irreducible and generic,
\item the quotient $V/\soc(V)$ has no generic subquotients (or, equivalently, the top derivative
$(V/\soc(V))^{(n)}$ vanishes), and
\item $V$ is the sum of its finite length submodules.
\end{enumerate}
\end{definition}

We will need the following ``dual version'' of this (c.f.~\cite{emerton-helm}, Lemma 6.3.5).

\begin{lemma}
Let $\kappa$ be a $W(k)$-algebra that is a field, and let $V$ be a finite length
admissible $\kappa[G]$-module such that the cosocle of $V$ is absolutely irreducible and generic,
and such that $V^{(n)}$ is a one-dimensional $\kappa$-vector space.  Then the smooth $\kappa$-dual
of $V$ is essentially AIG.
\end{lemma}

It follows easily from the definitions (see~\cite{emerton-helm}, Lemma 3.2.3 for details) that
the only endomorphisms of an essentially AIG $\kappa[G]$-module $V$ are scalars.  In particular
such a $V$ is indecomposable and lies in a single Bernstein component.  If $V$ lies in
$\Rep_{W(k)}(G)_{[L,\pi]}$, then $A_{[L,\pi]}$ acts on $V$, and this action factors through a map
$f_V: A_{[L,\pi]} \rightarrow \kappa$. 

\section{Projective objects and their derivatives} \label{sec:projective}

We now recall some facts from~\cite{bernstein1} about certain projective objects
of $\Rep_{W(k)}(G)$ and the action of the Bernstein center on them.  Let
$(K,\tau)$ be a maximal distingushed cuspidal $k$-type of $G$ in the sense of~\cite{vig98},
IV.3.1B.

From such a type we can construct a projective $W(k)[G]$-module $\CP_{K,\tau}$
by a construction detailed in section 4 of~\cite{bernstein1} (see particularly Lemmas 4.7
and 4.8, and the paragraph immediately following them).  It will be useful later to note
that the construction of $\CP_{K,\tau}$ shows that $\CP_{K,\tau}$ admits
a natural surjection
$$\CP_{K,\tau} \rightarrow \cInd_K^G \tau.$$
 
Moreover, we have:

\begin{theorem} Let $\pi$ be an irreducible cuspidal $k$-representation of $G$ containing
the type $(K,\tau)$, and let $(L,\pi')$ be the supercuspidal support of $\pi$.
Then $\CP_{K,\tau}$ lies in the block $\Rep_{W(k)}(G)_{[L,\pi']}$.  Moreover, the
map 
$$A_{[L,\pi']} \rightarrow \End_{W(k)[G]}(\CP_{K,\tau})$$
is an isomorphism, and makes $\CP_{K,\tau}$ into an admissible $A_{[L,\pi]}$-module.
\end{theorem}
\begin{proof}
This is a composite of results from~\cite{bernstein1}, in particular Corollary 10.9
and Corollary 11.11.  Admissibility follows from Theorem 8.8.
\end{proof}

Let $E_{K,\tau}$ be the ring $\End_{W(k)[G]}(\CP_{K,\tau})$.  We identify
$A_{[L,\pi']}$ with $E_{K,\tau}$ for the remainder of this section.
Fix a prime ${\mathfrak p}$ of $E_{K,\tau}$, of residue field $\kappa({\mathfrak p})$
and consider the admissible $\kappa({\mathfrak p})[G]$-module
$\CP_{K,\tau} \otimes_{E_{K,\tau}} \kappa({\mathfrak p})$.  It is contained in
a finite collection of blocks of $\Rep_{\kappa({\mathfrak p})}(G)$ and is
therefore of finite length.

\begin{lemma}
The cosocle of $\CP_{K,\tau} \otimes_{E_{K,\tau}} \kappa({\mathfrak p})$ is absolutely
irreducible.
\end{lemma}
\begin{proof}
Let $\tilde C$ be the cosocle of $\CP_{K,\tau} \otimes_{E_{K,\tau}} \kappa({\mathfrak p})$.
It suffices to show that $\End_{\kappa({\mathfrak p})}(\tilde C)$ is isomorphic to
$\kappa({\mathfrak p})$.  Let $C$ be the image of $\CP_{K,\tau}$ in $\tilde C$.  Becuase
$\CP_{K,\tau}$ is projective, the surjection of $\CP_{K,\tau}$ onto $C$ gives a surjection
of $E_{K,\tau}$ onto $\End_{W(k)[G]}(C)$.  On the other hand, the annihilator of $C$ in $E_{K,\tau}$
is ${\mathfrak p}$, so the map $E_{K,\tau}/{\mathfrak p} \rightarrow \End_{W(k)[G]}(C)$
is an isomorphism.  Tensoring with $\kappa({\mathfrak p})$ completes the proof.
\end{proof}

\begin{proposition} \label{prop:AIG point}
The following are equivalent:
\begin{enumerate}
\item $\CP_{K,\tau} \otimes_{E_{K,\tau}} \kappa({\mathfrak p})$ has an absolutely irreducible
generic quotient.
\item The smooth $\kappa({\mathfrak p})$-dual of $\CP_{K,\tau} \otimes_{E_{K,\tau}} \kappa({\mathfrak p})$
is essentially AIG.
\end{enumerate}
\end{proposition}
\begin{proof}
It is clear that (2) implies (1).  Suppose conversely that (1) holds.  It suffices to show
that $\CP_{K,\tau} \otimes_{E_{K,\tau}} \kappa({\mathfrak p})$ has no irreducible generic 
subquotients other than its cosocle.  Let $\tilde C$ denote this cosocle, and suppose we
have:
$$M \subseteq \ker \left[\CP_{K,\tau} \otimes_{E_{K,\tau}} \kappa({\mathfrak p}) \rightarrow \tilde C\right],$$
such that $M$ has an irreducible generic quotient $\tilde C'$.  The action of $A_{[L,\pi']}$
on both $\tilde C$ and $\tilde C'$ is via the same map
$$A_{[L,\pi']} \cong E_{K,\tau} \rightarrow \kappa({\mathfrak p}).$$
Thus $\tilde C$ and $\tilde C'$ have the same supercuspidal support.  (This follows by
Theorem 12.2 of~\cite{bernstein1} if ${\mathfrak p}$ has characteristic $\ell$, and by
classical results of Bernstein-Deligne~\cite{BD} if ${\mathfrak p}$ has characteristic zero.)
As there is a unique isomorphism class of irreducible generic representation with given supercuspidal
support, $\tilde C$ and $\tilde C'$ are isomorphic.  The surjection of $M$ onto $\tilde C'$
thus gives a surjection of $M$ onto $\tilde C$.  The surjection
$$\CP_{K,\tau} \otimes_{E_{K,\tau}} \kappa({\mathfrak p}) \rightarrow \tilde C$$
then lifts (by projectivity of $\CP_{K,\tau}$) to a map 
$$\CP_{K,\tau} \otimes_{E_{K,\tau}} \kappa({\mathfrak p}) \rightarrow M.$$
Composing this with the inclusion of $M$ in $\CP_{K,\tau} \otimes_{E_{K,\tau}} \kappa({\mathfrak p})$
gives a nonzero endomorphism of $\CP_{K,\tau} \otimes_{E_{K,\tau}} \kappa({\mathfrak p})$
that is zero on the cosocle $\tilde C$.  This is impossible since we have
$$\End_{\kappa({\mathfrak p})}(\CP_{K,\tau} \otimes_{E_{K,\tau}} \kappa({\mathfrak p}))
\cong E_{K,\tau} \otimes_{E_{K,\tau}} \kappa({\mathfrak p}) = \kappa({\mathfrak p}).$$
\end{proof}

This motivates the following definition:
\begin{definition} A point ${\mathfrak p}$ of $\Spec E_{K,\tau}$ is called
an {\em essentially AIG point} if $\CP_{K,\tau} \otimes_{E_{K,\tau}} \kappa({\mathfrak p})$
satisfies the equivalent conditions of Proposition~\ref{prop:AIG point}.
\end{definition}

Let ${\mathfrak m}$ and ${\mathfrak p}$ be points of $\Spec E_{K,\tau}$,
such that ${\mathfrak m}$ is in the closure of ${\mathfrak p}$.  
Then there is a discrete valuation ring $\OO \subset \kappa({\mathfrak p})$
that dominates ${\mathfrak m}$; that is, such that $\OO$ contains the image
of $E_{K,\tau}$ in $\kappa({\mathfrak p)}$, and the preimage of the maximal
ideal $\tilde {\mathfrak m}$
of $\OO$ under the map $E_{K,\tau} \rightarrow \OO$ is equal to
${\mathfrak m}$.

Let ${\tilde C}$ be the unique irreducible generic representation of $G$
over $\kappa({\mathfrak p})$ on which $E_{K,\tau}$ acts via the natural map
$E_{K,\tau} \rightarrow \kappa({\mathfrak p})$.  Then the supercuspidal
support of ${\tilde C}$ is $\OO$-integral, and so ${\tilde C}$ is
$\OO$-integral as well.  

Consider the smooth $\kappa({\mathfrak p})$-dual ${\tilde C}^{\vee}$
of $\tilde C$.
Then by~\cite{emerton-helm}, Proposition 3.3.2,
there is an $\OO$-lattice $C^{\vee}$ in ${\tilde C}^{\vee}$ such that
the reduction $C^{\vee} \otimes_{\OO} \OO/{\tilde {\mathfrak m}}$ has
an absolutely irreducible generic socle.  If we let $C$ be the smooth $\OO$-dual
of $C^{\vee}$, then $C \otimes_{\OO} \OO/{\tilde {\mathfrak m}}$ is an $\OO$-lattice
in $\tilde C$ with an absolutely irreducible generic quotient.  Denote this quotient
by $\overline{C}$.

\begin{lemma} \label{lemma:generize}
If ${\mathfrak m}$ is an essentially AIG point, then so is ${\mathfrak p}$.
\end{lemma}
\begin{proof}
Note that $E_{K,\tau}$ acts on $\overline{C}$ via the map $E_{K,\tau} \rightarrow
\kappa({\mathfrak m})$, and so $\overline{C}$ is the unique absolutely irreducible
generic representation of $G$ over $k$ on which $E_{K,\tau}$ acts via this map.
Since ${\mathfrak m}$ is an essentially AIG point, we have a map $\CP_{K,\tau} \rightarrow \overline{C}$.
Projectivity of $\CP_{K,\tau}$ lets us lift this to a map $\CP_{K,\tau} \rightarrow C$, and
hence to a map $\CP_{K,\tau} \rightarrow {\tilde C}$.  Thus ${\mathfrak p}$ is an essentially
AIG point.
\end{proof}

A partial converse also holds:

\begin{lemma} \label{lemma:specialize}
If ${\mathfrak p}$ is an essentially AIG point, and $C \otimes_{\OO} \OO/{\tilde {\mathfrak m}}$
is absolutely irreducible, then ${\mathfrak m}$ is an essentially AIG point.
\end{lemma}
\begin{proof}
If ${\mathfrak p}$ is an essentially AIG point, then we have a map
$\CP_{K,\tau} \rightarrow {\tilde C}$.  The image of $\CP_{K,\tau} \otimes_{E_{K,\tau}} \OO$
in $\tilde C$ is an admissible $\OO[G]$-submodule of $\tilde C$ and thus a $G$-stable
$\OO$-lattice in $\tilde C$.  As $C/{\tilde {\mathfrak m}}C$ is absolutely irreducible,
such an $\OO$-lattice must be homothetic to $C$.  We thus get a surjection
of $\CP_{K,\tau} \otimes_{E_{K,\tau}} \OO$ onto $C$.  Composing with the map
$C \rightarrow \overline{C}$ yields a nonzero map $\CP_{K,\tau} \rightarrow \overline{C}$
with kernel ${\mathfrak m}$, and the result follows.
\end{proof}

Of course, neither lemma is of much use without knowing some essentially AIG points
of $\Spec E_{K,\tau}$.  The following result provides such points:

\begin{lemma} \label{lemma:cusp AIG}
Let ${\mathfrak p}$ be a point of $\Spec E_{K,\tau}$, and suppose that there
exists an irreducible cuspidal representation $\pi$ of $G$ over $\kappa({\mathfrak p})$
on which $E_{K,\tau}$ acts via the natural map $E_{K,\tau} \rightarrow \kappa({\mathfrak p})$.
Then ${\mathfrak p}$ is an essentially AIG point.
\end{lemma}
\begin{proof}
If ${\mathfrak p}$ has characteristic zero, then $\pi$ is supercuspidal and is therefore
the unique irreducible representation with supercuspidal support $(G,\pi)$, and hence the
only irreducible representation on which $E_{K,\tau}$ acts via $E_{K,\tau} \rightarrow \kappa({\mathfrak p})$.
It follows that $\CP_{K,\tau} \otimes_{E_{K,\tau}} \kappa({\mathfrak p})$ is isomorphic
to $\pi$ in this case.

We may thus assume that ${\mathfrak p}$ has characteristic $\ell$.  Then $\pi$ is a cuspidal
$k$-representation of $G$ with supercuspidal support $(L,\pi')$.  It follows that $\pi$
contains the cuspidal type $(K,\tau)$.  In particular we have a map 
$$\cInd_K^G \tau \rightarrow \pi.$$

On the other hand, composing with the natural surjection of
$\CP_{K,\tau}$ onto $\cInd_K^G \tau$ yields a nonzero map
of $\CP_{K,\tau}$ onto $\pi$, and this map is annihilated by ${\mathfrak p}$.  
Tensoring with
$\tkappa({\mathfrak p})$ gives a surjection of $\CP_{K,\tau} \otimes_{E_{K,\tau}} \tkappa({\mathfrak p})$
onto $\pi$, proving the claim.
\end{proof}

Combining the above observations, we may now prove:

\begin{proposition}
Every point ${\mathfrak m}$ of $\Spec E_{K,\tau}$ is an essentially AIG point.
\end{proposition}
\begin{proof}
We will prove this by constructing points
${\mathfrak p}$ and ${\mathfrak m'}$ of $\Spec E_{K,\tau}$, such that
\item ${\mathfrak m}$ and ${\mathfrak m'}$ are in the closure of ${\mathfrak p}$,
and the unique irreducible generic representation on which $E_{K,\tau}$ acts
via $E_{K,\tau} \rightarrow \kappa({\mathfrak m'})$ is cuspidal. 
Then ${\mathfrak m'}$ is essentially AIG by Lemma~\ref{lemma:cusp AIG}, and
so ${\mathfrak p}$ is an essentially AIG point by Lemma~\ref{lemma:generize}.

Thus, if we can further arrange that the pair $({\mathfrak m},{\mathfrak p})$
satisfies the hypotheses of Lemma~\ref{lemma:specialize}, the claim will follow.
To do this we make use of the recent classification of modular representations of
$G$ due to M{\'{\i}}nguez-S{\'e}cherre~\cite{M-S}.  Let $\overline{C}$ be the unique absolutely
irreducible generic epresentation on which $E_{K,\tau}$ acts via
$E_{K_{\tau}} \rightarrow \kappa({\mathfrak m})$.

By~\cite{M-S}, Th{\'e}or{\`e}me 9.10, the representation $\overline{C}$ is
the parabolic induction of the tensor product
$$L(\Delta_1) \otimes \dots \otimes L(\Delta_r),$$
where each $\Delta_i$ is a {\em segment} in the sense of~\cite{M-S},
D{\'e}finition 7.1, and $L(\Delta_i)$ is the representation attached to $\Delta_i$
by~\cite{M-S}, D{\'e}finition 7.5.  Moreover, the segments $\Delta_i$ and $\Delta_j$
are not linked (\cite{M-S}, D{\'e}finition 7.3) for any $i,j$, and the multiset
$\{\Delta_1, \dots, \Delta_r\}$ is {\em aperiodic} in the sense of~\cite{M-S}, D{\'e}finition 9.7.

Let $\tkappa$ be the field of fractions of $\kappa({\mathfrak m})[T_1, \dots, T_r]$, and
let $\chi_1, \dots, \chi_r$ be the unramified characters $F^{\times} \rightarrow \tkappa^{\times}$
such that $\chi_i$ takes the value $T_i$ on a uniformizer $\unif_F$ of $F$.  For each $i$,
let $\Delta'_i$ be the twist $\Delta_i \otimes \chi_i$, and let $\tilde C$
be the parabolic induction of the tensor product
$$L(\Delta'_1) \otimes \dots \otimes L(\Delta'_r).$$
The segments $\Delta'_i$ and $\Delta'_j$ are not linked for any $i,j$, so $\tilde C$
is absolutely irreducible and generic.  Let $\OO$ be a valuation ring of $\tkappa$
that dominates the maximal ideal $\<T_i - 1\>$ in $\kappa({\mathfrak m})[T_1, \dots, T_r]$.  Then $\tilde C$
contains a stable $\OO$-lattice $C$, and the reduction $C \otimes_{\OO} \kappa({\mathfrak m})$
is isomorphic to $\overline{C}$.  In particular $\tilde C$ lies in the same block as $\overline{C}$
and thus admits an action of $E_{K,\tau}$; let ${\mathfrak p}$ be the kernel of this action.  It
is clear that $\kappa({\mathfrak p}) = \tkappa$.  Our construction shows that ${\mathfrak m}$
and ${\mathfrak p}$ satisfy the hypotheses of Lemma~\ref{lemma:specialize}, so that ${\mathfrak m}$
is an essentially AIG point if ${\mathfrak p}$ is.

It remains to construct an essentially AIG point ${\mathfrak m}'$ in the closure of ${\mathfrak p}$.
We first recall that, in the notation of~\cite{M-S}, section 6.2,
any cuspidal representation of $\GL_n(F)$ has the form
$\St(\pi_0,m)$ for some integer $m$ and supercuspidal representation $\pi_0$ of $\GL_{\frac{n}{m}}(F)$.
Since there is a cuspidal representation in $\Rep_k(G)_{[L,\pi]}$, we must have $\pi$
intertially equivalent to $\pi_0^{\otimes m}$ for a suitable $\pi_0$ and $m$.  On the other hand, because
$\overline{C}$ lies in $\Rep_k(G)_{[L,\pi]}$, the segments $\Delta_i$ must have the form
$[a_i,b_i]_{\pi_i}$ for all $i$, where $\pi_i$ is a cuspidal representation.  Then $\pi_i$ has
the form $\St(\pi_0,m_i) \otimes \chi'_i$ for some integer $m_i$ and some unramified character $\chi'_i$.

The supercuspidal support of $\St(\pi_0,m_i)$, considered
as a multiset of supercuspidal representations, is given by 
$$\pi_0, \pi_0 \otimes \nu, \dots, \pi_0 \otimes \nu^{m_i - 1},$$
where $\nu$ is the unramified character attached to $\pi_0$ by~\cite{M-S}, section 5.2.
It follows that the supercuspidal support of $L(\Delta_i)$ is given by
the union of the multisets
$$\pi_0 \otimes \chi'_i \otimes \nu^j, \pi_0 \otimes \chi'_i \otimes \nu^{j+1}, \dots, \pi_0 \otimes \chi'_i \otimes \nu^{j+ m_i - 1}$$
as $j$ ranges from $a_i$ to $b_i$.  Thus, if we choose suitable unramified characters
$\chi''_i: F^{\times} \rightarrow \kappa({\mathfrak m})^{\times}$,
and set $\Delta''_i = \Delta_i \otimes \chi''_i$
we can arrange that the parabolic induction of the tensor product
$$L(\Delta''_1) \otimes \dots \otimes L(\Delta''_r)$$
has the same supercuspidal support as the cuspidal representation $\St(\pi_0,m)$.

Observe that the natural map $E_{K,\tau} \rightarrow \kappa({\mathfrak p})$ factors
through the inclusion:
$$\kappa({\mathfrak m})[T_1^{\pm 1}, \dots, T_r^{\pm 1}] \rightarrow \kappa{(\mathfrak p)}.$$
Let $c_i = \chi''_i(\unif_F)$, and let ${\mathfrak m}'$ be the preimage of the maximal ideal
$\<T_i - c_i\>$ of $\kappa({\mathfrak m})[T_1^{\pm 1}, \dots, T_r^{\pm 1})$ in $E_{K,\tau}$.
Then ${\mathfrak m}'$ is the kernel of the map giving the action of $E_{K,\tau}$ on $\St_{\pi_0,m}$.
In particular ${\mathfrak m}'$ is an essentially AIG point in the closure of ${\mathfrak p}$,
and the result follows.
\end{proof}

\begin{corollary}
The space $\CP_{K,\tau}^{(n)}$ is locally free of rank one over $E_{K,\tau}$.
\end{corollary}
\begin{proof}
Admissibility of $\CP_{K,\tau}$ as an $E_{K,\tau}[G]$-module (and hence as an $A_{[L,\pi]}[G]$-module)
implies by \cite{emerton-helm}, 3.1.14, that $\CP_{K,\tau}^{(n)}$ is finitely 
generated over $E_{K,\tau}$.
On the other hand, for any prime ${\mathfrak p}$ of $E_{K,\tau}$ we have
$$\CP_{K,\tau}^{(n)} \otimes_{E_{K,\tau}} \kappa({\mathfrak p}) \cong
(\CP_{K,\tau} \otimes_{E_{K,\tau}} \kappa({\mathfrak p}))^{(n)}.$$
As ${\mathfrak p}$ is an essentially AIG point the right-hand side is a one-dimensional
$\kappa({\mathfrak p})$ vector space; the result follows.
\end{proof}

\section{Admissibility of $W_{L,\pi}$} \label{sec:admissible}

Now let $L$ be an arbitrary Levi subgroup of $\GL_n$, and let $\pi$
be a supercuspidal representation of $L$.  Our goal is to use the results of
the previous section to show that the module $W_{L,\pi}$ is admissible over $A_{L,\pi}$.
We must first relate $A_{L,\pi}$ to the spaces $\CP_{K,\tau}$ studied above;
we do so by invoking results of~\cite{bernstein1} which we now recall:

Section 11 of~\cite{bernstein1}
attaches to the pair $(L,\pi)$ a Levi subgroup $M^{\max}$ of $G$, containing $L$
and an irreducible cuspidal $k$-representation $\pi^{\max}$ of $M^{\max}$, that is maximal
(with respect to a certain partial order defined in~\cite{bernstein1}) among
pairs $(M,\pi')$ of a Levi $M$ of $G$ and a cuspidal representation $\pi'$ of $M$ over $k$
such that $(M,\pi')$ has inertial supercuspidal support $(L,\pi)$.)  The Levi $M^{\max}$
is isomorphic to a product
$$\GL_{n_1}(F) \times \dots \times \GL_{n_r}(F),$$
and the representation $\pi^{\max}$ is a tensor product of irreducible cuspidal representations
$\pi_i$ of $\GL_{n_i}(F)$ over $k$.

Each $\pi_i$ contains a maximal distinguished $k$-type $(K_i,\tau_i)$.  We can thus form
the representation $\otimes_i \CP_{K_i,\tau_i}$ of $M^{\max}$.  Parabolically inducing
from $M^{\max}$ to $G$ then yields a projective module
$$\CP_{M^{\max},\pi^{\max}} := i_P^G \otimes_i \CP_{K_i,\tau_i}$$
for a suitable parabolic $P$ of $G$ with Levi $M^{\max}$.

Let $E_i$ be the endomorphism ring $E_{K_i,\tau_i}$.  Then the tensor product
$\otimes_i E_i$ acts on $\otimes_i \CP_{K_i,\tau_i}$ and thus, by functoriality of parabolic
induction, on $\CP_{M^{\max},\pi^{\max}}$.  Moreover, we have:

\begin{proposition} \label{prop:saturation}
The action of $A_{[L,\pi]}$ on $P_{M^{\max},\pi^{\max}}$ factors uniquely through the map
$$\otimes_i E_i \rightarrow \End_{W(k)[G]}(\CP_{M^{\max},\pi^{\max}}),$$
and the resulting map $A_{[L,\pi]} \rightarrow \otimes_i E_i$ is an injection.
Moreover, if $x$ is an element of $\otimes_i E_i$ such that $\ell^r x$ is in the image
of $A_{[L,\pi]}$, then $x$ is also in the image of $A_{[L,\pi]}$.
\end{proposition}
\begin{proof}
This is immediate from~\cite{bernstein1}, Corollary 11.11 and Proposition 11.5.
\end{proof}

We are now in a position to show:

\begin{theorem} \label{thm:whittaker endomorphisms}
The natural map $A_{[L,\pi]} \rightarrow \End_{W(k)[G]}(W_{[L,\pi]})$ is an isomorphism.
\end{theorem}
\begin{proof}
We already know that this map becomes an isomorphism after tensoring with $\overline{\CK}$,
and hence after inverting $\ell$.  In particular (as $A_{[L,\pi]}$ and $W_{[L,\pi]}$ have no
$\ell$-torsion), the map in question is injective.  It thus suffices to show surjectivity.

We have an isomorphism:
$$\CP_{M^{\max},\pi^{\max}}^{(n)} \cong \Hom_{W(k)[G]}(W_{[L,\pi]}, \CP_{M^{\max},\pi^{\max}}),$$
and therefore an action of $\End_{W(k)[G]}(W_{[L,\pi]})$ on $\CP_{M^{\max},\pi^{\max}}^{(n)}$.
On the other hand, the multiplicativity of the top derivative with respect to parabolic induction
yields an isomorphism of $\CP_{M^{\max},\pi^{\max}}^{(n)}$ with the tensor product
of the invertible $E_i$-modules $\CP_{K_i,\tau_i}^{(n_i)}$.  The action of
$\End_{W(k)[G]}(W_{[L,\pi]})$ on $\CP_{M^{\max},\pi^{\max}}^{(n)}$ thus yields a map
$$\End_{W(k)[G]}(W_{[L,\pi]}) \rightarrow \otimes_i E_i$$
that extends the inclusion of $A_{L,\pi}$ into $\otimes_i E_i$.  Such an extension is
necessarily injective, as $\otimes_i E_i$ is $\ell$-torsion free.

We thus have inclusions:
$$A_{[L,\pi]} \subseteq \End_{W(k)[G]}(W_{[L,\pi]}) \subseteq \otimes_i E_i.$$
Now let $x$ be an element of $\End_{W(k)[G]}(W_{[L,\pi]})$.  Then for some positive integer $a$,
$\ell^a x$ lies in $A_{[L,\pi]}$.  But then by Proposition~\ref{prop:saturation}, 
$x$ must have been an element of $A_{[L,\pi]}$ to start with.  The result follows.
\end{proof}

\begin{proposition} \label{prop:whittaker admissibility}
The module $W_{[L,\pi]}$ is admissible as an $A_{[L,\pi]}[G]$-module.
\end{proposition}
\begin{proof}
Theorem~\ref{thm:whittaker endomorphisms}, together with Lemma~\ref{lemma:freeness},
show that $W_{[L,\pi]}^{(n)}$ is free of rank one over $A_{[L,\pi]}$.  Let $x$ be
an element of $W_{[L,\pi]}^{(n)}$ that generates $W_{[L,\pi]}^{(n)}$ as an $A_{[L,\pi]}$-module,
and let $\tx$ be an element of $W_{[L,\pi]}$ that maps to $x$ via the natural surjection:
$$W_{[L,\pi]} \rightarrow W_{[L,\pi]}^{(n)}.$$

Let $W'$ be the $W(k)[G]$-submodule of $W_{[L,\pi]}$ generated by $\tx$.  The inclusion
of $W'$ in $W_{[L,\pi]}$ gives an inclusion of $(W')^{(n)}$ in $W_{[L,\pi]}^{(n)}$
whose image contains $x$; it follows that $(W')^{(n)}$ is equal to $W_{[L,\pi]}^{(n)}$.
Thus $(W_{[L,\pi]}/W')^{(n)} = 0$.  But the latter is naturally isomorphic to
$\Hom_{W(k)[G]}(W_{[L,\pi]},W_{[L,\pi]}/W')$, so we must have $W' = W_{[L,\pi]}$.
In particular $W_{[L,\pi]}$ is finitely generated over $W(k)[G]$ and is thus
(by~\cite{bernstein1}, Corollary 12.4)
admissible over $A_{[L,\pi]}$.
\end{proof}

We now have the following ``structure theory'' for essentially AIG $k[G]$-modules (or, more
precisely, their duals):
\begin{proposition} 
Let $\kappa$ be a $W(k)$-algebra that is a field, and let
$f: A_{[L,\pi]} \rightarrow \kappa$ be a map of $W(k)$-algebras.
Then $[W_{[L,\pi]} \otimes_{A_{[L,\pi]},f} \kappa]^{\vee}$ is essentially AIG,
where the superscript $\vee$ denotes $\kappa$-dual.

Conversely, let $V$ be an object of
$\Rep_{W(k)}(G)_{[L,\pi]}$ that is the smooth $\kappa$-dual of an
essentially AIG $\kappa[G]$-module.  Then $V$ is 
a quotient of $W_{[L,\pi]} \otimes_{A_{[L,\pi]},f} \kappa$ for some 
$f: A_{[L,\pi]} \rightarrow \kappa$.
\end{proposition}
\begin{proof}
By Theorem~\ref{thm:whittaker endomorphisms}, together with the fact that
$W_{[L,\pi]}^{(n)}$ is free of rank one over $\End_{W(k)[G]}(W_{[L,\pi]})$,
we find that $[W_{[L,\pi]} \otimes_{A_{[L,\pi]},f} \kappa]^{(n)}$ is a one-dimensional
$\kappa$-vector space.  Moreover, if $\pi$ is any nonzero quotient of
$W_{[L,\pi]} \otimes_{A_{[L,\pi]},f} \kappa$, then $\Hom_{W(k)[G]}(W,\pi)$ is nonzero,
so that $\pi^{(n)}$ is nonzero.  Thus, by~\cite{emerton-helm}, Lemma 6.3.5, the $\kappa$-dual of 
$W_{[L,\pi]} \otimes_{A_{[L,\pi]},f} \kappa$
is essentially AIG.

Now let $V$ be a $\kappa[G]$ module in $\Rep_{W(k)}(G)_{[L,\pi]}$, such that $V$ is the smooth dual of an
essentially AIG module.  Then $V$ has an absolutely irreducible generic cosocle $V_0$, and
$A_{[L,\pi]}$ acts on $V_0$ by a map $f: A_{[L,\pi]} \rightarrow \kappa$.
As the only endomorphisms of $V$ are scalars, $A_{[L,\pi]}$ acts on $V$ via $f$, as well.

As $V_0$ is irreducible and generic, and is an object of $\Rep_{W(k)}(G)_{[L,\pi]}$,
we have a nonzero map of $W_{[L,\pi]}$ onto $V_0$; projectivity of $W_{[L,\pi]}$
lifts this to a map $W_{[L,\pi]} \rightarrow V$.  Let ${\mathfrak p}$ be the kernel of $f$;
then ${\mathfrak p}$ acts by zero on
$V$, so the map $W_{[L,\pi]} \rightarrow V$ descends to a nonzero map:
$$W_{[L,\pi]}/{\mathfrak p} W_{[L,\pi]} \rightarrow V;$$
as $V$ is a $\kappa$-vector space, and $A/{\mathfrak p}$ is contained in $\kappa$, this extends to
a map:
$$W_{[L,\pi]} \otimes_A \kappa \rightarrow V$$
whose composition with the map $V \rightarrow V_0$ is nonzero.
This map is necessarily surjective
(if not, its image would be contained in the kernel of the map $V \rightarrow V_0$,
as $V_0$ is the cosocle of $V$). 
\end{proof}

\begin{corollary}
Every essentially AIG $\kappa[G]$-module has finite length.
\end{corollary}

This verifies a conjecture of~\cite{emerton-helm} (see the remarks after Lemma 3.2.8).

\section{co-Whittaker modules} \label{sec:whittaker}

We now describe a sense in which $W_{[L,\pi]}$ can be thought of as a ``universal object''
over $A_{[L,\pi]}$.  This will turn out to be crucial to our reinterpretation of the
results of~\cite{emerton-helm}.

Let $A$ be a Noetherian $W(k)$-algebra; then any smooth $A[G]$-module carries an action of
the center of $\Rep_{W(k)}(G)$, and hence admits a Bernstein decomposition.  Thus
$\Rep_A(G)$ decomposes as a product of Bernstein components $\Rep_A(G)_{[L,\pi]}$.

\begin{definition} An object $V$ of $\Rep_A(G)_{[L,\pi]}$ is a {\em co-Whittaker
A[G]-module} if the following hold:
\begin{enumerate}
\item $V$ is admissible as an $A[G]$-module.
\item $V^{(n)}$ is free of rank one over $A$.
\item If ${\mathfrak p}$ is a prime ideal of $A$, with residue field
$\kappa({\mathfrak p})$, then the $\kappa({\mathfrak p})$-dual
of $V \otimes_A \kappa({\mathfrak p})$ is essentially AIG.
\end{enumerate}
If $V$ and $V'$ are co-Whittaker $A[G]$-modules, we say that $V$ dominates $V'$
if there exists a surjection from $V$ to $V'$.
\end{definition}

One motivation for this definition is that the families of admissible representations
that correspond to Galois representations in the sense of section 6.2 of~\cite{emerton-helm} are
co-Whittaker modules.  We will discuss this further in the section~\ref{sec:langlands}.

\begin{proposition}
Let $V$ be a co-Whittaker $A[G]$-module.  Then the natural map
$$A \rightarrow \End_{A[G]}(V)$$
is an isomorphism.
\end{proposition}
\begin{proof}
By localizing at each prime ideal of $A$, it suffices to consider the case in which $A$ is a local
ring.  Lemma 6.3.2 of~\cite{emerton-helm} then shows that $V$ is generated by a certain
submodule ${\mathfrak J}(V)$ that is stable under any endomorphism of $V$.  The proof then
proceeds exactly as in the proof of part (3) of Proposition 6.3.4 of~\cite{emerton-helm}.
\end{proof}

If $V$ is a co-Whittaker $A[G]$-module of $\Rep_A(G)_{[L,\pi]}$, then it admits
an action of $A_{[L,\pi]}$, which must arise from a unique map $f_V: A_{[L,\pi]} \rightarrow A$.

\begin{theorem} \label{thm:universal}
Let $A$ be a Noetherian $A_{[L,\pi]}$-algebra.
Then $W_{[L,\pi]} \otimes_{A_{[L,\pi]}} A$ is a co-Whittaker
$A[G]$-module.  Conversely, if $V$ is a co-Whittaker $A[G]$-module in $\Rep_A(G)_{[L,\pi]}$, then
$A$ is an $A_{[L,\pi]}$-algebra via $f_V: A_{[L,\pi]} \rightarrow A$, and
$W_{[L,\pi]} \otimes_{A_{[L,\pi]}} A$ dominates $V$. 
\end{theorem}
\begin{proof}
The module $W_{[L,\pi]}$ is admissible over $A_{[L,\pi]}[G]$, so 
$W_{[L,\pi]} \otimes_{A_{[L,\pi]}} A$ is admissible over $A$.
As $W_{[L,\pi]}^{(n)}$ is free of rank one over $A_{[L,\pi]}$, and the derivative operator
commutes with base change, $[W_{[L,\pi]} \otimes_{A_{[L,\pi]}} A]^{(n)}$ is free of
rank one over $A$.  Finally, as $W_{[L,\pi]} \otimes_{A_{[L,\pi]}} \kappa({\mathfrak p})$
is dual to an essentially AIG representation
for any prime ideal ${\mathfrak p}$ of $A_{[L,\pi]}$,
it follows that $W_{[L,\pi]} \otimes_{A_{[L,\pi]}} A$ is a co-Whittaker $A[G]$-module.

For the converse, choose a generator of $V^{(n)}$ as an $A$-module; such a generator
corresponds to a map $W_{[L,\pi]} \rightarrow V$.  This map induces a map
$W \otimes_{A_{[L,\pi]}} A \rightarrow V$ of $A[G]$-modules which we must prove is surjective.
Let $V'$ be the cokernel.  We have an exact sequence:
$$0 \rightarrow (W_{[L,\pi]} \otimes_{A_{[L,\pi]}} A)^{(n)} \rightarrow V^{(n)} \rightarrow (V')^{(n)}
\rightarrow 0$$
and the first horizontal map is surjective by construction.  Thus $(V')^{(n)} = 0$.
Assume $V'$ is nonzero; then
since $V'$ is admissible over $A[G]$ it has a quotient that is simple as an $A[G]$-module;
this is a non-generic quotient of $V \otimes_A \kappa({\mathfrak m)}$ for some maximal ideal
${\mathfrak m}$ of $A$.  But $V \otimes_A \kappa({\mathfrak m})$ has an absolutely irreducible
generic cosocle so this cannot happen.
\end{proof}

We conclude with a technical lemma which will be useful in section~\ref{sec:langlands}.

\begin{lemma} \label{lemma:interpolation}
Let $A$ be a reduced Noetherian $A_{[L,\pi]}$-algebra, and suppose that for each
minimal prime ${\mathfrak a}$ of $A$, we specify a quotient $V_{\mathfrak a}$
of $W_{[L,\pi]} \otimes_{A_{[L,\pi]}} \kappa({\mathfrak a})$.  Then there exists
a co-Whittaker $A[G]$-module $V$, unique up to isomorphism, such that
$V$ is $A$-torsion free, and
$V \otimes_A \kappa({\mathfrak a})$ is isomorphic to $V_{\mathfrak a}$
for all ${\mathfrak a}$.
\end{lemma}
\begin{proof}
Let $V$ be the image of the diagonal map:
$$W_{[L,\pi]} \otimes_{A_{[L,\pi]}} A \rightarrow \prod_{\mathfrak a} V_{\mathfrak a}.$$
The $A[G]$-module $V$ is clearly $A$-torsion free by construction,
and $V \otimes_A \kappa({\mathfrak a})$ is clearly isomorphic to $V_{\mathfrak a}$.
It thus suffices to show that $V$ is a co-Whittaker module.  But $V$ is a quotient of
the co-Whittaker module $W_{[L,\pi]} \otimes_{A_{L,\pi}} A$, so it suffices to
show that $V^{(n)}$ is free of rank one over $A$.  Certainly $V^{(n)}$ is cyclic,
so it suffices to show that $V^{(n)} \otimes_A \kappa({\mathfrak a})$ is nonzero
for all minimal primes $\kappa({\mathfrak a})$.  This is clear because
$V \otimes_A \kappa({\mathfrak a})$ is isomorphic to the essentially AIG
representation $V_{\mathfrak a}$.

As for the uniqueness claim, any such $V$ is dominated by
$W_{[L,\pi]} \otimes_{A_{[L,\pi]}} A$, and embeds in
$\prod_{\mathfrak a} V_{\mathfrak a}$.  Up to an automorphism of
$\prod_{\mathfrak a} V_{\mathfrak a}$, the composition of
$W_{[L,\pi]} \otimes_{A_{[L,\pi]}} A \rightarrow V$ with 
the embedding of $V$ in $\prod_{\mathfrak a} V_{\mathfrak a}$
is equal to the diagonal map, and thus identifies $V$ with the image of
this diagonal map.
\end{proof} 

\section{The local Langlands correspondence in families} \label{sec:langlands}

We now apply our results to the local Langlands correspondence in families
of~\cite{emerton-helm}, which we now recall:

\begin{theorem}[\cite{emerton-helm}, Theorem 6.2.1] \label{thm:E-H}
Let $A$ be a reduced complete Noetherian local $\ell$-torsion free
$W(k)$-algebra, with residue field $k$,
and let $\rho: G_F \rightarrow \GL_n(A)$ be a Galois representation.
Then there is, up to isomorphism, at most one admissible $A[G]$-module
$\pi(\rho)$ such that:
\begin{enumerate}
\item $\pi(\rho)$ is $A$-torsion free,
\item $\pi(\rho)$ is a co-Whittaker $A[G]$-module, and
\item for each minimal prime ${\mathfrak a}$ of $A$, the representation
$\pi(\rho)_{{\mathfrak a}}$ is $\kappa({\mathfrak a})$-dual to the representation
that corresponds to $\rho^{\vee}_{\mathfrak a}$ via the Breuil-Schneider
generic local Langlands correspondence.  (For a discussion of this
correspondence, see~\cite{emerton-helm}, section 4.2.)
\end{enumerate}
\end{theorem}

It is conjectured (\cite{emerton-helm}, Conjecture 1.1.3) that such a $\pi(\rho)$ exists
for any $\rho$.  The goal of this section is to reformulate this conjecture as a
relationship between the Bernstein center $A_{[L,\pi]}$ and the deformation theory of
Galois representations.

To motivate our reformulation, it will be useful to invoke the following result,
which will be proved in the forthcoming work~\cite{bernstein2}.  Note that
we only use this result for motivation, and that our main result (Theorem~\ref{thm:main})
does not depend on it.

\begin{proposition} \label{prop:conditional}
Fix a representation $\rhobar: G_F \rightarrow \GL_n(k)$,
and let $(R_{\rhobar}^{\Box},\rho^{\Box})$ be the universal framed deformation of $\rhobar$.
The ring $R_{\rhobar}^{\Box}$ is reduced and $\ell$-torsion free.
\end{proposition}

Granting this, the representation $\rho^{\Box}$ fits into the framework of Theorem~\ref{thm:E-H}.
Let us suppose that the family $\pi(\rho^{\Box})$ exists.  (If this is true, then
the argument of~\cite{emerton-helm}, Proposition 6.2.10, allows us to construct the family $\pi(\rho)$ 
for {\em any} deformation $\rho$ of $\rhobar$, essentially via ``base change''.)  
Let $\pibar$ be an irreducible representation of $G$ over $k$ whose supercuspidal
support corresponds to $\rhobar$ under the mod $\ell$ semisimple local Langlands
correspondence of Vign{\'e}ras~\cite{vigss}, and
choose a pair $(L,\pi)$ such that $\pibar$ lies in the block
$\Rep_{W(k)}(G)_{[L,\pi]}$.  Then every subquotient of $\pi(\rho^{\Box}) \otimes_{R_{\rhobar}^{\Box}} k$
has the same supercuspidal support as $\pibar$.

By~\cite{emerton-helm}, Theorem 6.2.1, the endomorphism ring $\End_{W(k)[G]}(\pi(\rho^{\Box}))$
is isomorphic to $R_{\rhobar}^{\Box}$.  As this ring is local (and thus has no nontrivial idempotents)
the action of the Bernstein center of $\pi(\rho^{\Box})$ must factor through $A_{[L,\pi]}$.
We thus obtain a map:

$$A_{[L,\pi]} \rightarrow R_{\rhobar}^{\Box}.$$

Note that (as $R_{\rhobar}^{\Box}$ is complete and local by definition), this map factors
through the completion $(A_{[L,\pi]})_{\mathfrak m}$, where ${\mathfrak m}$ is the
kernel of the action of $A_{[L,\pi]}$ on $\pibar$; this yields
a map:
$$\LL: (A_{[L,\pi]})_{\mathfrak m} \rightarrow R_{\rhobar}^{\Box}.$$

We can easily describe the effect of this map on the characteristic zero points of
$\Spec R_{\rhobar}^{\Box}$.  Let $\kappa$ be a complete field of characteristic zero that
contains $W(k)$.  A representation $\rho: G_F \rightarrow \kappa$,
determines, via (Tate normalized) local Langlands,
an irreducible representation $\Pi$ in $\Rep_{\kappa}(G)$.  If $\Pi$ lies in
$\Rep_{\kappa}(G)_{[L,\pi]}$, we let $f_{\rho}$ denote the map $A_{[L,\pi]} \rightarrow \kappa$
giving the action of $A_{[L,\pi]}$ on $\Pi$.  Note that this map depends only on the supercuspidal
support of $\Pi$, so that $f_{\rho}$ depends only on the semisimplification
of $\rho$.

Now let $\rho$ be a deformation of $\rhobar$ over $\kappa$.  A suitable choice of basis for $\rho$ 
determines a map $x: R_{\rhobar}^{\Box} \rightarrow \kappa$, such that
$\rho^{\Box}_x = \rho$.  It follows from~\cite{emerton-helm}, Theorem 6.2.6 that
every subquotient of $\pi(\rho^{\Box})_x$ then has supercuspidal support corresponding to $\rho$
under (Tate normalized) local Langlands.  On the other hand, $A_{[L,\pi]}$ acts on
$\pi(\rho^{\Box})_x$ via the composition
$$A_{[L,\pi]} \rightarrow (A_{[L,\pi]})_{\mathfrak m} \stackrel{\LL}{\rightarrow} R_{\rhobar}^{\Box}
\stackrel{x}{\rightarrow} \kappa.$$
It follows that this composition is equal to $f_{\rho}$.  Summarizing, we have:

\begin{proposition} \label{prop:interpolate}
Suppose that $\pi(\rho^{\Box})$ exists, and let $\kappa$ be a complete field of characteristic
zero containing $W(k)$.  Then the map:
$$\LL: (\Spec R_{\rhobar}^{\Box})(\kappa) \rightarrow (\Spec (A_{[L,\pi]})_{\mathfrak m})(\kappa)$$
takes the point $x$ in $(\Spec R_{\rhobar}^{\Box})(\kappa)$ that corresponds to a framed deformation $\rho$
to the $\kappa$-point $f_{\rho}$ of $\Spec (A_{[L,\pi]})_{\mathfrak m}$.
\end{proposition}

In less formal language, we will say that a map 
$(A_{[L,\pi]})_{\mathfrak m} \rightarrow R_{\rhobar}^{\Box}$ 
``interpolates the characteristic zero semisimple local Langlands
correspondence'' over $\Spec R_{\rhobar}^{\Box}$ if it satisfies the conclusion of Proposition~\ref{prop:interpolate}.

The following then follow easily from Proposition~\ref{prop:conditional},
together with the fact that $A_{[L,\pi]}$ is reduced and $\ell$-torsion free:

\begin{itemize}
\item There is at most one map:
$(A_{[L,\pi]})_{\mathfrak m} \rightarrow R_{\rhobar}^{\Box}$ that interpolates the characteristic zero
semisimple local Langlands correspondence; i.e. $\LL$ is determined by Proposition~\ref{prop:interpolate}
if it exists.
\item If $\rhobar$ is semisimple, then the map $\LL$ is injective.
\item The image of $\LL$ in $R_{\rhobar}^{\Box}$ is contained in the set of functions that are
``invariant under change of frame''.  More precisely, if $G^{\Box}$ is the formal completion of $\GL_n/W(k)$
at the $k$-point corresponding to the identity identity, then $G^{\Box}$ acts on $R_{\rhobar}^{\Box}$
by changing the frame, and the image of $\LL$ is contained in the invariants
$(R_{\rhobar}^{\Box})^{G^{\Box}}$ of this action.
\end{itemize}

This motivates the following conjecture, which (if we grant Proposition~\ref{prop:conditional})
is a consequence of the existence of $\pi(\rho^{\Box})$.

\begin{conj} \label{conj:weak}
For any $\rhobar$, there is a map:
$$\LL: (A_{[L,\pi]})_{\mathfrak m} \rightarrow R_{\rhobar}^{\Box}$$
that interpolates the characteristic zero semisimple local Langlands correspondence.
\end{conj}

It is also tempting, given the naturality of the rings $A_{[L,\pi]}$ and $(R_{\rhobar}^{\Box})^{G^{\Box}}$,
to formulate a stronger version of this conjecture when $\rhobar$ is semisimple:

\begin{conj} \label{conj:strong}
If $\rhobar$ is semisimple, there is a unique isomorphism:
$$\LL: (A_{[L,\pi]})_{\mathfrak m} \rightarrow (R_{\rhobar}^{\Box})^{G^{\Box}}$$
that interpolates the characteristic zero semisimple local Langlands correspondence.
\end{conj}

\begin{remark} \label{rem:interpolation}
Conjecture~\ref{conj:strong} can be thought of as a geometric interpolation
of the semisimple local Langlands correspondence in characteristic zero.  Indeed,
the points of $(A_{[L,\pi]})_{\mathfrak m}$ correspond to supercuspidal
supports of representations in $A_{[L,\pi]}$.  On the other hand, let
$\kappa$ be a complete field of characteristic zero containing $W(k)$, and let
$x$ and $y$ be two $\kappa$-points of $\Spec R_{\rhobar}^{\Box}$.  It follows from the theory
of pseudocharacters that $f(x) = f(y)$ for all $G^{\Box}$-invariant
elements $f$ of $R_{\rhobar}^{\Box}$ 
if, and only if the Galois representations $\rho_x$ and $\rho_y$
have the same semisimplification.  One can thus regard the $\kappa$-points of
$\Spec (R_{\rhobar}^{\Box})^{G^{\Box}}$ as parameterizing ``lifts of $\rhobar$
to $\kappa$ up to semisimplification''; from this point of view the conjectured
bijection:
$$(\Spec (R_{\rhobar}^{\Box})^{G^{\Box}})(\kappa) \rightarrow (\Spec (A_{[L,\pi]})_{\mathfrak m})(\kappa)$$
is simply the map that takes a Galois representation (up to semisimplification) to
the corresponding supercuspidal support.
\end{remark}

\begin{remark}
Conjecture~\ref{conj:strong} can also be thought of as an integrality result
for certain Bernstein center elements defined by Galois-theoretic considerations.
For instance, Chenevier has shown that for any $\tau \in W_F$, there is a unique element
$f_{\tau}$ of the center of $\Rep_{\overline{\CK}}(G)$ whose action on irreducible admissible
representations $\Pi$ of $G$ over $\overline{\CK}$ is given in the following way:
let $\rho$ be the Frobenius-semisimple representation of $W_F$ over $\overline{\CK}$ corresponding to $\Pi$;
then $f_{\tau}$ acts on $\Pi$ by the trace of $\rho(\tau)$.  (We refer the reader to~\cite{chenevier}, Proposition 3.11
for details.)

Conjecture~\ref{conj:strong} would imply that Bernstein center elements defined in this way
actually live in the completion of the center of $\Rep_{W(k)}(G)$ at {\em every} maximal
ideal ${\mathfrak m}$, and thus lie in the center of $\Rep_{W(k)}(G)$.
\end{remark}

The existence of $\pi(\rho^{\Box})$, together with Proposition~\ref{prop:conditional},
would imply Conjecture~\ref{conj:weak} for $\rhobar$.

Our main result is that the converse also holds; that is, that 
Conjecture~\ref{conj:weak} {\em implies}
to the existence of $\pi(\rho)$ for all deformations $\rho$ of $\rhobar$ over
suitable coefficient rings.  

\begin{theorem} \label{thm:main}
Fix an $A$ and a $\rho$ as in Theorem~\ref{thm:E-H},
and let $\rhobar = \rho \otimes_A k$.  Suppose that Conjecture~\ref{conj:weak}
holds for $\rhobar$.
Then $\pi(\rho)$ exists.
\end{theorem}
\begin{proof}
A choice of basis for $\rho$ yields a map $R_{\rhobar}^{\Box} \rightarrow A$.
Composition with the map $\LL$ makes $A$ into an $A_{[L,\pi]}$-algebra.
By Lemma~\ref{lemma:interpolation}, to construct $\pi(\rho)$,
it suffices to show that for all ${\mathfrak a}$,
the representation $\pi_{\eta}$ that
is $\kappa({\mathfrak a})$-dual to the representation attached to
$\rho^{\vee}_{\mathfrak a}$ by the Breuil-Schneider correspondence is a quotient
of $W_{[L,\pi]} \otimes_{A_{[L,\pi]}} \kappa({\mathfrak \eta})$.   It suffices
to show that $\pi_{\eta}$ is dual to an essentially AIG representation, and that
its supercuspidal support corresponds to that of $W_{[L,\pi]} \otimes_{A_{[L,\pi]}}
\kappa({\mathfrak \eta}).$  The representation $\pi_{\eta}$ is dual
to an essentially AIG representation because any representation
arising from the Breuil-Schneider correspondence is essentially AIG;
see~\cite{emerton-helm}, Corollary 4.3.3.  It has the correct supercuspidal support
because the map $\LL$ interpolates the semisimple local Langlands correspondence.
\end{proof}

\begin{remark} In~\cite{bernstein2}, we will show that both Conjecture~\ref{conj:weak}
and Conjecture~\ref{conj:strong} hold after inverting $\ell$, and also hold if $\ell$
is a {\em banal} prime; that is, if $1, q, \dots, q^n$ are distinct modulo $\ell$.  We
will also establish Conjecture~\ref{conj:weak} when $n=2$ and $\ell$ is odd.  (This last
result relies on forthcoming work of Paige~\cite{paige} on the structure of projective
envelopes of representations of $\GL_n(\FF_q)$.)

In particular the results of~\cite{bernstein2}, together with Theorem~\ref{thm:main},
will establish the existence of $\pi(\rho)$ for any two-dimensional representation
$\rho$ of $G_F$ over a complete local ring $A$ of odd residue characteristic that
satisfies the conditions of Theorem~\ref{thm:E-H}.
\end{remark}

\begin{example} \label{ex:banal}
To illustrate the power of Theorem~\ref{thm:main}, we sketch a proof
of Conjecture~\ref{conj:strong} (and hence of the existence of a local Langlands correspondence
in families), in the case when $\ell$ is banal.  (The details of this argument will appear
in~\cite{bernstein2}.)

The basic ideas are as follows:
When $\ell$ is banal, the pair $(M^{\max},\pi^{\max})$ described in section~\ref{sec:admissible}
is equal to $(L,\pi)$.  It follows that in this case the algebra $A_{[L,\pi]}$ is equal
to the subalgebra $C_{[L,\pi]}$ of $A_{[L,\pi]}$ described explicitly in section 11 of~\cite{bernstein1}.
Explicitly, if we fix supercuspidal representations $\pi_1, \dots, \pi_s$ that are in distinct inertial
equivalence classes, such that $\pi$ is inertially equivalent to the tensor product
$$\pi_1^{r_1} \otimes \dots \otimes \pi_s^{r_s},$$
and for each $i$ we fix a lift $\tpi_i$ of $\pi_i$ to a representation over $\CK$,
then $A_{[L,\pi]}$ is freely generated by elements $\Theta_{i,j}$, where $i$ ranges from
one to $s$ and $j$ ranges from one to $r_i$, together with the inverses of $\Theta_{i,r_i}$ for each $i$.

An irreducible $\Pi$ in $\Rep_{\overline{\CK}}(G)_{[L,\pi]}$ in this case has supercuspidal support
inertially equivalent to the tensor product
$$\tpi_1^{r_1} \otimes \dots \otimes \tpi_s^{r_s},$$
or more precisely of the form
$$\{\tpi_i \otimes \chi_{ij} : 1 \leq i \leq s; 1 \leq j \leq r_i\}$$
for unramified characteris $\chi_{ij}$.  If we set $\chi_{ij}(\unif_F) = \alpha_{ij}$,
then an element $\Theta_{i,k}$ acts on $\Pi$ via the $k$-th elementary symmetric function on the scalars
$\alpha_{i,j}$ for $i$ fixed and $1 \leq j \leq r_i$.

On the Galois side, we have for each $i$ a $\trho_i: G_F \rightarrow \GL_{n_i}(W(k))$ corresponding to 
$\tpi_i$ via local Langlands.  An adaptation of an argument due to Clozel-Harris-Taylor shows that any 
deformation of the semisimple $\rhobar$ corresponding to $\pi$ has the form:
$$\bigoplus_i V_i \otimes \trho_i,$$
where $V_i$ is an $r_i$-dimensional representation of $G_F$ such that the action of the inertia group $I_F$ on $V_i$ factors
through its pro-$\ell$ quotient.  The fact that $\ell$ is banal then implies that the action of $I_F$
on $V_i$ is unipotent; more precisely the action of $G_F$ on $V_i$ can be given by a pair of
$r_i$ by $r_i$ matrices $F_i, \sigma_i$ with $F_i \sigma_i F_i^{-1} = \sigma_i^q$, and $\sigma_i$ unipotent.
From this one easily
deduces that $(R_{\rhobar}^{\Box})^{G^{\Box}}$ is freely topologically generated by the coefficients
of the characteristic polynomial of $F_i$ for each $i$, together with the inverses of the $\det F_i$.

Moreover, it is easy to see that the map that takes $\Theta_{i,k}$ to $(-1)^k$ times the coefficient of $X^{r_i - k}$
in the characteristic polynomial of $F_i$ interpolates the characteristic zero local Langlands correspondence,
and thus provides the desired isomorphism:
$$(A_{[L,\pi]})_{\mathfrak m} \rightarrow (R_{\rhobar}^{\Box})^{G^{\Box}}.$$

In constrast to this elementary computation, providing a direct construction of $\pi(\rho^{\Box})$
by elementary methods seems to be very hard.  In particular it is easy to construct $\pi(\rho^{\Box})$
over each irreducible component of $\Spec R_{\rhobar}^{\Box}$, but ``gluing'' these branches of the family together
in the unique way that satisfies condition (2) of Theorem~\ref{thm:E-H} requires verifying the existence of
suitable compatibilities between these branches over the intersections of the irreducible components.  In general
these intersections can be quite singular, and the compatibilities become increasingly difficult to verify as the
number of components grows, rendering this approach essentially hopeless for large $n$.
\end{example}


\begin{thebibliography}{}
\bibitem[BD]{BD}
J. Bernstein and P. Deligne, \emph{Le ``centre'' de Bernstein}, in
\emph{Representations des groups redutifs sur un corps local, Traveaux en cours},
(P. Deligne ed.), Hermann, Paris, 1--32.

\bibitem[BH]{BH-whittaker}
C. Bushnell and G. Henniart, \emph{Generalized Whittaker models and the Bernstein center,}
Amer. J. Math. {\bf 125} (2003), no. 3, 513--547.

\bibitem[BK]{BK}
C. Bushnell and P. Kutzko, \emph{The admissible dual of $\GL(N)$ via compact open subgroups,}
Annals of Math. Studies {\bf 129}, Princeton University Press, Princeton, 1993.

\bibitem[C]{chenevier}
G. Chenevier, \emph{Une application des vari{\'e}ti{\'e}s de Hecke des groups unitaires,}
preprint, 2011.

\bibitem[D]{dat-integral}
J.-F. Dat, \emph{Integral elements in Bernstein's center}, preprint, 2006.

\bibitem[E]{emerton-lg}
M. Emerton, \emph{Local-global compatibility in the $p$-adic Langlands programme for $\GL_2/\QQ$,}
preprint, 2011.

\bibitem[EH]{emerton-helm}
M. Emerton and D. Helm, \emph{The local Langlands correspondence for $\GL_n$ in families},
preprint, 2011, {\bf arXiv:1104.0321}.

\bibitem[H1]{bernstein1}
D. Helm, \emph{The Bernstein center of the category of smooth $W(k)[\GL_n(F)]$-modules},
preprint, 2012, {\bf arxiv:1201.1874}.

\bibitem[H2]{bernstein2}
D. Helm, \emph{Towards a Galois-theoretic interpretation of the integral Bernstein center},
in preparation.

\bibitem[JPS]{JPS}
H. Jacquet, I. I. Piatetski-Shapiro, and J. Shalika, \emph{Conducteur des repr{\'e}sentations du
groupe lineare}, Math. Annalen {\bf 256} (1981), no. 2, 199--214.

\bibitem[MS]{M-S}
A. M{\'{\i}}nguez and V. S{\'e}cherre, \emph{Repr{\'e}sentations lisses modulo $\ell$ de $\GL(m,D)$},
preprint, 2011, {\bf arxiv:1110.1467}.

\bibitem[P]{paige}
D. Paige, {\emph The projective envelope of a cuspidal representation of a finite linear group},
in preparation.

\bibitem[V1]{vig98}
M.-F. Vign{\'e}ras, \emph{Induced $R$-representations of $p$-adic reductive groups},
Selecta Math. {\bf 4} (1998), no. 4, 549--623.

\bibitem[V2]{vigss}
M.-F. Vign{\'e}ras, \emph{Correspondance de Langlands semi-simple pour $\GL(n,F)$ modulo $\ell \neq p$},
Invent. Math. {\bf 144} (2001), 177-223.



\end{thebibliography}
\end{document}